\documentclass[10pt,a4paper]{article}
\usepackage{amssymb}
\newtheorem{thm}{Theorem}
\newtheorem{defn}{Definition}
\newtheorem{cor}{Corollary}
\newtheorem{lem}{Lemma}
\newtheorem{prop}{Proposition}
\newtheorem{ex}{Example}
\newtheorem{rk}{Remark}

\usepackage{amssymb}

\def \C{{\! \rm \
I \!\!\!C}}

\def \R {{\! \rm \ I \!R}}

\def \N {{\! \rm \ I \!N}}

\def \Z {{\! \rm Z\! \!Z}}

\def\otherterm#1{{\it#1}}
\def \e {{\epsilon}}
\newcommand{\tr} {\hbox{tr}}

\def \Ci {{C^\infty}}
\def \l {{\lambda}}
\def \tr {{\rm tr}}
\def \Cl {{C\ell}}

\def\sres{{\rm sres}}
\def \endsquare{ $\sqcup \!\!\!\! \sqcap$ }

\begin{document} 
\title{\bf The  logarithmic residue density of a generalised Laplacian}
\author{  Jouko Mickelsson and Sylvie Paycha  }
\maketitle
\section*{Abstract}
We show that the residue density of the logarithm of a generalised Laplacian
on a closed manifold 
defines an invariant polynomial valued differential form. We express it in 
terms of a
finite sum of    residues  of classical pseudodifferential symbols. In  the
case of the square of a Dirac operator, these formulae provide a 
pedestrian proof of the   
Atiyah-Singer formula for a pure Dirac operator in dimension $4$ and for a
twisted Dirac operator on  a flat space of any dimension. These correspond to 
special cases of a more general formula  by S. Scott and D. Zagier announced 
in \cite{Sc2} and  to appear in \cite{Sc3}. In our approach, which is of perturbative
nature, we  use  either
a Campbell-Hausdorff formula derived by Okikolu or a non commutative Taylor
type 
 formula. 
\section*{Introduction}
The noncommutative residue on classical pseudodifferential operators has a
notable property, locality i.e., it corresponds to a residue density ${\rm
  res}_x(A)\, dx$ integrated over an $n$-dimensional (closed) manifold $M$:
\begin{equation}\label{eq:resintro} {\rm res}(A)= \int_M {\rm res}_x(A)\, dx; \quad  {\rm
  res}_x(A)\, dx:= \frac{1}{(2\pi)^n}\int_{\vert \xi\vert=1}{\rm tr}\left( \sigma_{-n}(A)(x,
\xi)\right)\, d\xi,\end{equation} where $\sigma(A)(x, \xi)$ stands
for the local symbol of $A$, $\sigma_{-n}(A)(x, \xi)$ for the $-n$-th
homogeneous part of its symbol,  with $(x,\xi)$ varying in the cotangent bundle of $M$ and
where tr stands for the fibrewise trace. 
As it was observed in \cite{O2}, this  extends to logarithms $A=\log Q$ of an elliptic
pseudodifferential operator $Q$  of positive order with appropriate spectral cut (we call it
admissible). \\  Exponentiating ${\rm res}(\log (Q))$ leads to the  residue determinant ${\rm det}_{\rm res}(Q):= e^{{\rm
  res}(\log Q)}$ first
introduced by Wodzicki in the case of zero order operators (see
e.g. the survey \cite{Ka}) and further   extended by Scott \cite{Sc1}  to  elliptic
pseudodifferential operator $Q$  with appropriate spectral cut of positive
order. Further logarithmic structures were since then investigated in \cite{Sc2}
in relation with topological quantum field theory.\\
Here, we show that  the logarithmic residue density for a generalised Laplacian
$Q$,  $${\rm res}_x(\log Q)\, dx:=\frac{1}{(2\pi)^n}\,\left[ \int_{S_x^*M}{\rm
  tr}\left(\sigma_{-n}(\log Q)(x, \xi) \right)\, d\xi\right]\, dx, $$  which
defines  an invariant polynomial valued  form in the
sense of Weyl (Theorem \ref{thm:resinvpol}). It then follows from Gilkey's invariance theory, that it can be
expressed in terms of Pontryagin and Chern classes.\\
The presence of a logarithm makes the actual computation of a  logarithmic residue
density difficult. However, observing    that  the symbol of a generalised
Laplacian reads
$\sigma(Q)= \vert  \xi\vert^2+\sigma_{<2}(Q)$ where $\sigma_{<2}(Q)$ is of
order smaller than $2$ enables us  to carry out computations by  means of 
 a noncommutative Taylor type
formula (Theorem \ref{thm:NCTaylorloggsymbol}) or a   
Campbell-Hausdorff formula  (Theorem \ref{thm:CHlogresdensity}), both of
which provide ways to compare $\sigma(\log Q)$   with 
$\log\left( \vert  \xi\vert^2\right) $.\\ A similar procedure applies on the
operator level to compute the integrated logarithmic density $\int_M {\rm
res}_x\left(
\sigma(\log Q)\right)\, dx$  in the
special case of the square of a  twisted Dirac operator $D_W$ acting on a
 twisted $\Z_2$-graded spinor
bundle $E=S\otimes W$. Indeed, combining  the Lichnerowicz
formula (\ref{eq:Bochner})  
(which compares $D_W^2$ with a Laplace-Beltrami  operator
$\Delta^E=\left(\nabla^E\right)^*\nabla^E$ built from the underlying
 connection $\nabla^E$ on $E$)   
 with
a Campbell-Hausdorff formula (which compares $\log D_W^2$ with 
$\log \Delta^E$)
 yields an expression for  the integrated 
logarithmic   (super-) residue density  ${\rm sres}(\log D_W^2)$ in terms  of ${\rm
sres}(\log (\Delta^E))$ and a finite number of (super-) residues of
classical operators involving the curvature of $\nabla^E$ 
(Theorem \ref{thm:indtwisted}). This integrated 
logarithmic   (super-) residue density turns out to be 
proportional to the index of the chiral Dirac operator $D_W^+$ \footnote{This was observed independently by
 S. Scott in \cite{Sc1} and the
  second author in some unpublished lecture notes delivered in Colombia.}
  (Theorem \ref{thm:superresidue}):\begin{equation}\label{eq:indintro}
  {\rm ind}(D^+)=-\frac{1}
{2}\int_M {\rm sres}_x(\log D_W^2)\, dx.\end{equation}  
 Thus, locality in the  Atiyah-Singer index theorem is  closely related to
the local property of the  noncommutative residue. \\ 
We  compute the index   in two concrete examples, first for a twisted Dirac
operator on a flat space (Theorem \ref{thm:flattorus}) along the lines
described above using a Campbell-Hausdorff formula, then for a pure
Dirac operator in dimension four using a Taylor type formula. For the second
example we first derive simple formulae 
(see Proposition \ref{prop:srespartialGamma}) for (super-) residues of
certain expressions involving the derivatives of the Christoffel symbols,
which can then be used to derive the index in dimension four. 
We recover this
way, the Atiyah-Singer index theorem for a pure Dirac operator on 
a four dimensional spin manifold.\\
 With the perturbative approach adopted here using either a Campbell-Hausdorff or
 a noncommutative Taylor formula, we were unfortunately  unable to derive the general
 Atiyah-Singer formula   announced in \cite{Sc2}
 and to appear \cite{Sc3}. This perturbative approach
 nevertheless provides a  pedestrian proof in the cases investigated
 here and useful  intermediate results such as Theorem \ref{thm:resinvpol}
 and Theorem \ref{thm:NCTaylorloggsymbol}, which we find are 
  of
  interest in their own right.
  
\section*{Notations}
Given  an even $n= 2p$-dimensional real  oriented euclidean
vector space $V,$  there is a unique $\Z_2$-graded complex Clifford module $S=S^+\oplus
S^-$,
the spinor module, such that the complex Clifford algebra
 $C(V)\otimes \C\simeq {\rm End}(S)$ and  ${\rm dim}(S)=
2^p$. An auxillary linear complex space $W$ yields a $\Z_2$-graded twisted Clifford module
$E=S\otimes W$. 
Let  \begin{eqnarray*}
c:
\Lambda V&\to& C(V)\\ 
e_{i_1}\wedge \cdots \wedge e_{ik}&\mapsto & c(e_{i_1}) \cdots c( e_{i_k} )
\end{eqnarray*}
be the  quantisation map. To simplify notations we set
$\gamma_j= c(e_j)$, so that the grading operator reads
 $\Gamma=i^p \, \gamma_1\cdots \gamma_n$. Notice that $\Gamma^2= Id$.  The
 cyclicity of the trace combined with the Clifford relations imply that the
  supertrace ${\rm str}:= {\rm
   tr}\circ \Gamma$  on ${\rm End}(E)$  satisfies the following property for a
   matrix $M\in {\rm End}(W)$ viewed as an element of ${\rm End}(E)$:
\begin{equation}\label{eq:strgamma}  
     {\rm str}(M\,\gamma_{i_1}\cdots \gamma_{i_k})= 0  \quad {\rm if} \quad k <
n,
\quad {\rm str}(M\,\gamma_{1}\cdots \gamma_{n})=  
(-2i)^p {\rm tr}(M),\\
\end{equation}
since ${\rm dim}{\rm End}({\cal S})=2^p$.
On the other hand, setting 
 \begin{equation}\label{eq:sigmaij}\sigma_{ij}= 
 \frac{1}{8}[\gamma_i, \gamma_j]= \frac{1}{4}
  \gamma_i\gamma_j\quad{\rm if}\quad i\neq j,\end{equation} we have that for any
permutation $\tau\in \Sigma_n$ with signature $\vert \tau\vert$:
\begin{equation}\label{eq:strprodsigma}{\rm str}\left(\sigma_{\tau(1)\tau(2)} \,\sigma_{\tau(3)\tau(4)}\cdots
  \sigma_{\tau(n-1)\tau(n)}\right)= \frac{(-1)^{\vert \tau\vert}}{4^p}\, {\rm
  str}\left(\gamma_1\cdots \gamma_n\right)= \frac{(-1)^{\vert \tau\vert}\,
   (-i)^p}{2^p}.\end{equation} 
 These constructions carry out to bundles, for which we abusively use the same
 notations. Let  $E= S\otimes W= E^+\oplus E^-$, with $E^+=S^+\otimes W$,
 $E^-=S^-\otimes W$ be a twisted $\Z_2$-graded spinor bundle
 over an even $n=2p$-dimensional  closed Riemannian manifold $M$, with
 auxillary bundle $W$  equipped with a connection $\nabla^W$. \\ 
 Let, for a vector bundle $F$ over $M$, 
 $\Cl(M, F)$ denote the algebra of classical
pseudodifferential operators acting on the space $\Ci(M, F)$ of 
 smooth sections of the vector bundle
$F$.  \\
 Let $D=\sum_{i=1}^n c(e_i)\, \nabla^S_{e_i}\in \Cl(M, S)$ be the
   Dirac operator, where $\nabla^S$ is the spinor
   connection, where  $c$ stands for the Clifford
   multiplication and $\{e_i, i=1,\cdots, n\}$ 
   for an orthonormal tangent frame on $M$.  In local coordinates we shall also
   write $\gamma_i$ for $c(e_i)$.\\
   Let  $\nabla^E:= \nabla^S\otimes 1+ 1\otimes \nabla^W$ be a 
  connection on the twisted bundle $E= S\otimes W$ and let  $D_W=\sum_{i=1}^n
  c(e_i)\, \nabla^E_{e_i}\in \Cl(M< E)$ be the corresponding twisted 
   Dirac operator.   The chiral Dirac operators $D_W^+$
   and its formal adjoint $D_W^-$ act from $\Ci(M, E^+)$ to $\Ci(M, E^-)$ and
   conversely. 
\section{The logarithmic residue density   as an invariant polynomial} 
 Following Gilkey's notations (see (2.4.3) in \cite{G}), for a multi-index $\alpha=(\alpha_1 \dots \alpha_s)$ we introduce formal variables 
$g_{ij/\alpha}=\partial_\alpha g_{ij} $ for the partial derivatives of the metric tensor $g$ on $M$ and the connection $\omega$ on the external bundle.
Let us set $${\rm ord}(g_{ij/\alpha})=\vert \alpha\vert = \alpha_1 +\dots + \alpha_s; \quad {\rm
  ord}(\omega_{i/\beta})=\vert \beta\vert. $$Inspired by Gilkey  (see (1.8.18)
and (1.8.19) in \cite{G}), we set the following definition. \begin{defn}
 We call a
classical   operator $A\in \Cl(M, E)$  of order $a$
{\bf geometric}, if in any local trivialisation, the
homogeneous components $\sigma_{a-j}(A)$ are
 homogeneous of order $j$
in the  jets of the metric and of  the connection. 
\end{defn}
\begin{rk}
A                 differential operator  $A=\sum_{\vert \alpha\vert \leq
    a} c_\alpha(x)\, \partial_x^\alpha \in \Cl(M, E)$ is geometric if $ c_\alpha(x) $
is homogeneous of order $j=a-\vert \alpha\vert$
in the  jets of the metric and of  the connection
$\nabla^W$. Here we use the standard notation ${\partial_x}^{\alpha} = \partial_{\alpha_1} \dots \partial_{\alpha_s}.$
\end{rk}
\begin{ex} 
The Laplace Beltrami operator 
$$\Delta_g=-\frac{1}{\sqrt g} \sum_{i=1, j=1}^n \partial_i
 \left( \sqrt g\, g^{ij}\partial_j\right)$$ has this property. \\
More generally, formula  (2.4.22) in \cite{G}  shows
that $\Delta_p=d_{p-1}\delta_{p-1}+\delta_p d_p$ on $p$-forms, where $\delta_k= (-1)^{nk+1}\, \star_{n-k} \, d_{n-k-1} \star_{k+1}$, is a geometric operator. Indeed, each
derivative applied to $\star$ reduces the order of differentiation by $1$ and increases the order in the jets of the
 metric by $1$.
\end{ex}
\begin{ex} The square of the twisted Dirac operator
$$D_W^2= -\sum_{ij} g^{ij} \left( \nabla^E_i \nabla^E_j+ \sum_{k}\Gamma_{ij}^k
  \nabla^E_k\right)+
\sum_{i<j} c(dx^i)c(dx^j) [\nabla^E_i, \nabla^E_j]$$
has this property.
\end{ex}
Geometric operators form an algebra. 
\begin{lem}\label{lem:prodgeometric} The product of two geometric
 operators $A$ and $ B$ in $\Cl(M, E)$
  is again a geometric operator.
\end{lem}
{\bf Proof:} Since the product $A\, B$ has symbol $$\sigma(A\, B)\sim\sum_{\alpha } \frac{(-i)^{|\alpha|}}{\alpha!}
\partial_\xi^\alpha\sigma(A)\, \partial_x^\alpha\sigma(B),$$
we have $$\sigma_{a+b-k}(AB)=\sum_{\vert \alpha\vert +i+j=k } \frac{(-i)^{|\alpha|}}{\alpha!}
\partial_\xi^\alpha\sigma_{a-i}(A)\,
\partial_x^\alpha\sigma_{b-j}(B)$$
where $a$ is the order of $A$, $b$ the order of $B$.
Thus, if $\sigma_{a-i}(A)$ and $\sigma_{b-j}(B)$ are homogeneous of degree $i$
and $j$ respectively in the jets of the metric and the
connection, $\sigma_{a+b-k}(AB)$ is homogeneous of degree $i+j+\vert
\alpha\vert=k$.\endsquare\\ \\ 
Following \cite{BGV}, we call generalised Laplacian on $E$ a second order
differential operator acting on $\Ci(M, E)$ with leading symbol $\vert
\xi\vert^2$. Since generalised Laplacians are expected to be geometric (see
the examples in the first section), we assume generalised Laplacians are
geometric
without further specification. Note that a generalised Laplacian is admissible
(see  the Appendix). The following result  provides a way to build
  families  of geometric
operators.\begin{prop}\label{prop:geometricgerms} Let $Q\in \Cl(M, E)$ be a
  generalised Laplacian with  spectral cut\footnote{See the Appendix.} $\theta$. \\
Then  for any geometric operator $A$ in $ \Cl(M, E)$, the family  $A(z):=A\,
Q_\theta^z$ is a  family of geometric
operators. 
\end{prop}
{\bf Proof:} By Lemma \ref{lem:prodgeometric}, it is sufficient to prove the
result for $A=I$. For convenience, we drop the explicit mention of the
spectral cut. \\ 
Since  $$Q^z=\frac{1}{2i\pi}\int_\Gamma \lambda^z\, (Q-\lambda)^{-1}\,
  d\lambda,$$
where $\Gamma$ is a contour described in the Appendix (see formula (\ref{eq:contourGamma})),  we need to investigate the resolvent $R(Q, \lambda)=
  (Q-\lambda)^{-1}$, the homogeneous components   $\sigma_{2-j}(R(Q, \lambda))$ of the
  symbol of which are defined inductively on $j$ by 
\begin{eqnarray}\label{eq:resolvent}\sigma_{-2}(R(Q, \lambda))&=&
  (\sigma_2(Q))^{-1}, \nonumber\\ \sigma_{-2-j}(R(Q,
  \lambda))&=&-\sigma_{-2}(
R(Q, \lambda))\sum_{k+l+\vert
  \alpha\vert=j,l<j}\frac{(- i)^{\vert \alpha\vert}}{\alpha!} \, D_\xi^\alpha
\sigma_{2-k}(Q) D_x^\alpha \sigma_{-2-l}(R(Q, \lambda)) \nonumber\\
\end{eqnarray}
Using (\ref{eq:resolvent}), one shows by induction on $j$  that 
$\sigma_{-2-j}(R(Q,  \lambda) )$ is a finite sum of expressions of the type
$$(-i)^{\vert \alpha\vert}\, (\vert\xi\vert^2-\lambda)^{-1-k}\,
D_\xi^{\alpha_1}D_x^{\beta_1}\sigma_{2-l_1}(Q)\cdots
D_\xi^{\alpha_k}D_x^{\beta_k}\sigma_{2- l_k}(Q), \quad  \vert l\vert +\vert
\alpha\vert=j, \quad \vert \alpha\vert=\vert \beta\vert.$$
Inserting this in   \begin{equation}\label{eq:complexpower}
 \sigma_{2\,z-j}(Q^z)(x,\xi)=- \frac{1}{2i\pi}\int_\Gamma \lambda^z \,
 \sigma_{-2-j}(R(Q, \lambda))(x, \xi)\, d\lambda, 
\end{equation} and applying repeated integrations
by parts to compute the Cauchy integrals:
$$ \frac{1}{2i\pi} \int_\Gamma \lambda^z\, (\vert\xi\vert^2-\lambda)^{-k-1}\,
d\lambda
=(-1)^{k}\,\frac{ z(z+1)\cdots (z+(k-1))}{k!  }\, \vert\xi\vert^{2\,(z-k)}, $$
combination of symbols  of the type
\begin{equation}\label{eq:symbQz} \vert\xi\vert^{q(z-k)}
D_\xi^{\alpha_1}D_x^{\beta_1}\sigma_{q-l_1}(Q)(x, \xi)\cdots
D_\xi^{\alpha_k}D_x^{\beta_k}\sigma_{q- l_k}(Q)(x, \xi), \quad  \vert l\vert +\vert
\alpha\vert=j, \quad \vert \alpha\vert=\vert \beta\vert.
\end{equation}
Since $\sigma_{2- l}(Q)$ is homogeneous of order $l$ in the jets of the metric
and the connection, it follows that  for any complex number $z$, the symbol
$\sigma_{2z-j}(Q^z)$ is homogeneous of order $j$ as  a linear combination of  products of  homogeneous expressions  of order $j_i$ in the jets of the metric
and the connection such that $j_1+\cdots +j_k=j$. 
\endsquare
\\ \\ The notion of geometric operator extends to logarithms of admissible
operators defined in  the Appendix. 
\begin{defn} We say that the
  logarithm $\log_\theta A$ (see  formula (\ref{eq:sigmaj0log}) in the Appendix) of an admissible 
   operator $A\in \Cl(M, E)$  of order $a$ with spectral cut $\theta$  is {\bf
    geometric} 
 if in any local trivialisation, the
homogeneous components $\sigma_{-j,0}(\log_\theta A)$ are
 homogeneous of order $j$
in the  jets of the metric and of  the connection. 
\end{defn}
\begin{rk} This can be generalised to any log-polyhomogeneous operator $A$ of
  order $a$ and  
  logarithmic degree $k$ by requiring
  that all the coefficients $\sigma_{a-j,l}(A), l\in \{0, \cdots,
  k\}$ in the logarithmic expansion of the symbol be homogeneous of order
  $a-j$ 
in the  jets of the metric and of  the connection. 
\end{rk}
\begin{cor} The logarithm
 of a generalised Laplacian  is a geometric operator. 
\end{cor}
{\bf Proof:} Again, we drop the explicit mention of the spectral cut. 
Differentiating  (\ref{eq:symbQz}) w.r.t. $z$
at zero shows that 
 $\sigma_{-j,0}(\log  Q)(x, \xi)$ is a linear
combination of symbols  of the type
$$\vert\xi\vert^{-2\,k}
D_\xi^{\alpha_1}D_x^{\beta_1}\sigma_{2-l_1}(Q)(x, \xi)\cdots
D_\xi^{\alpha_k}D_x^{\beta_k}\sigma_{2- l_k}(Q)(x, \xi), \quad  \vert l\vert +\vert
\alpha\vert=j, \quad \vert \alpha\vert=\vert \beta\vert.$$
Hence the symbol
$\sigma_{-j,0}(\log Q)$ is homogeneous of order $j$ as  a linear combination of  products of  homogeneous expressions  of order $j_i$ in the jets of the metric
and the connection such that $j_1+\cdots +j_k=j$.
\endsquare
\begin{rk} This is a particular instance of a more general result, namely that
  the derivative $A^\prime(0)$ at zero of a holomorphic germ $A(z)\in \Cl(M, E)$ of geometric
  operators  around zero,  is also geometric.
\end{rk}
Adopting Gilkey's  notations  (see \cite{G} par. 2.4) let us  denote by ${\cal P}_{n,k,p}^{g,
  \nabla^W}  
$   (which we write ${\cal P}_{n,k,p}^{g}  
$ if $ E=S$) 
    the linear space consisting of $p$- form valued invariant\footnote{By
  invariant we mean that they agree in any coordinate system around $x_0$ which is
  normalised w.r. to the point $x_0$, i.e. such that
  $g_{ij}(x_0)=\delta_{i-j}$ and $\partial_kg_{ij}(x_0)=0$.}  polynomials that are homogeneous of
order $k$ in the jets of the metric\footnote{The order in the jets of the
  metric is defined by ord$\left(\partial_x^\alpha g_{ij}\right)=
\vert\alpha\vert$.} and of  the connection
$\nabla^W$. 
\begin{ex} The scalar curvature $r_M$ belongs to $ {\cal P}_{n ,2,0}^{g}$
  since it
  reads \hfill \break \noindent $r_M=2\sum_{i, j}\left(\partial^2_{i,j}g_{ij}-
  \partial^2_{i,i}g_{jj}\right)$
in Riamann normal coordinate system. 
\end{ex}  
\begin{thm} \label{thm:resinvpol}
The logarithmic residue density of a  generalised Laplacian $Q$ on $E$ 
\begin{equation}\label{eq:logresdensity} R_n(x,Q):= {\rm res}_x(\log_\theta Q)\, dx:=\frac{1}{(2\pi)^n}\,\left[ \int_{S_x^*M}{\rm
  tr}\left(\sigma_{-n, 0}(\log_\theta Q)(x, \xi) \right)\, d\xi\right]\, dx
\end{equation}   is an invariant polynomial in ${\cal P}_{n,n,n}^{g, \nabla^W}$. \\
 It is  
 generated by Pontrjagin forms of the tangent bundle and Chern forms on the auxillary bundle.
\end{thm}
{\bf Proof:}
By Proposition \ref{prop:geometricgerms}  the  logarithm (we drop the spectral
cut) $\log Q$ is geometric, so that  $\sigma_{-n}(\log Q)(x, \xi)$ is homogeneous of
degree $n$ in the jets of the metric and the connection. Integrating this
expression in $\xi$ on the unit cosphere shows that  the residue
density 
lies in   $ {\cal P}_{n ,n,n}^{g, \nabla}$.\\
 Since $ {\cal P}_{n ,n,n}^{g,
  \nabla^W}$ is a polynomial in the $2$-jets of the metric and the one jets of
the auxillary connection, it is 
 generated by Pontrjagin forms of the tangent bundle (see Theorem 2.6.2 in
\cite{G}) and Chern forms on the auxillary bundle.
\endsquare
\begin{rk}The logarithmic residue density    is clearly additive on direct sums
  $E_1\oplus E_2$ of vector bundles over a closed manifold $M$
 $$R_n(x,Q_1\oplus Q_2)= R_n(x,Q_1)+ R_n(x,Q_2)$$
but there is  a priori no reason why it should be  multiplicative on tensor products
$E_1\otimes E_2\to M_1\times M_2$ of vector bundles  $E_i$ over  closed
manifolds $M_i$.
\end{rk}

 \section{The logarithmic residue density via
  the Campbell-Hausdorff formula}
  The Campbell-Hausdorff formula provides a first approach to compute a
  local logarithmic residue density. 
By the results of Okikiolu \cite{O1}, for two admissible classical
pseudodifferential operators with scalar leading symbols $A$ and $B$ in
$\Cl(M, E)$ and under suitable technical assumptions on their spectrum to
ensure that their logarithms are well-defined, 
we have \begin{equation}\label{eq:CHop}
\log(AB)\sim \log A+\log B+\sum_{k=2}^{\infty } C^{(k)}(\log A, \log
  B),
\end{equation}
where $C^{(k)}( \log A, \log B)$ are  Lie monomials given by:
\begin{equation} \label{eq:Ck} C^{(k)}(P,Q):=\sum_{j=1}^\infty \frac{(-1)^{j+1}}{(j+1)} \sum
\frac{({\rm Ad }_P)^{\alpha_1}({\rm Ad }_Q)^{\beta_1}\cdots
 ({\rm Ad }_P)^{\alpha_j}({\rm Ad }_Q)^{\beta_j}}{(1+ \sum_{l=1}^j
  \beta_l)\alpha_1!\cdots\alpha_j! \beta_1!\cdots\beta_j!}(Q),  \end{equation}
which vanish
if $\beta_j>1$ or if $\beta_j=0$ and $\alpha_j>1$ and  with the inner sum running over $j-$ tuples of pairs $(\alpha_i,
  \beta_i)$ such that $\alpha_i+\beta_i>0$ and $\sum_{i=1}^j
  \alpha_i+ \beta_i= k$. Here ${\rm Ad}_P(Q)=[P,Q]$ and the symbol $\sim$ means that for any integer $n$ the
  difference
    \begin{equation}\label{eq:Fn}
F_n(A, B):=\log(AB)-\log A-\log B-\sum_{k=2}^{n+1 } C^{(k)}(\log A, \log
  B)
\end{equation}  is of order smaller than $-n$. The fact that the leading symbols are scalar
    ensures that the order of $C^{(k)} (\log A, \log B)$ decreases as $k$
    increases and hence a good control on the asymptotics as a result of the
    fact that the 
     adjoint operations ${\rm ad}_{\log A}$ and ${\rm ad}_{\log B}$ decrease
     the order by one unit. 
 \begin{prop}  \label{prop:CH}\cite{O1} Let $A, B\in \Cl(M, E)$ be invertible elliptic
  operators with scalar leading symbols  such
  that $A, B$ and their product $A\, B$ have well defined logarithms, then 
$F_n(A, B)$ defined as in (\ref{eq:Fn}) is trace-class so that
its Wodzicki residue vanishes.  Both its trace and its residue vanish  (see also \cite{Sc1}), 
\begin{equation}\label{eq:resCH} {\rm res}\left(\log (AB)-\log A-\log B\right)= 0.
\end{equation}
\end{prop} 
The proof in \cite{O1}  is based on a similar
expansion on the        level of  symbols which we now describe for future
use.
We consider the algebra ${\cal F}S(U)$ of formal  symbols on an open subset $U$
of $\R^n$ equipped with the symbol product $\star$
$$\sigma_1\star \sigma_2(x, \xi) = \sum_{\alpha \in
\N^n}\frac{(-i)^{\vert
  \alpha\vert}}{\alpha!}\partial_{\xi}^{\alpha}\sigma_1(x,\xi)\partial_x^{\alpha}\sigma_2(x,\xi).$$
Let $\{\sigma, \tau\}_\star:= \sigma\star \tau-\tau \star \sigma$ denote the
associated star bracket.  
\begin{ex} It $\sigma $ is polynomial, this formal power series of symbols
  with decreasing order becomes a finite sum, as in the
  following example of interest to us: 
$$ \{ \vert \xi\vert^{2}, \tau\}_\star
=( L_x+\Delta_x)\tau 
$$
where we have set 
\begin{equation}\label{eq:LxDeltax} 
L_x:= -2i  \sum_{a=1}^n \xi_a \partial_{x_a}\quad {\rm  and }\quad \Delta_x:= -\sum_{a=1}^n
\partial^2_{x_a}. \end{equation} 
\end{ex}
We define  ${\rm ad}^{*k}_\sigma$ by induction on $k$ setting ${\rm
  ad}^{*0}_\sigma(\tau)= \tau$ and 
${\rm ad}^{*(k+1)}_{\sigma}(\tau):= \{\sigma, {\rm ad}^{*k}_{\sigma}(\tau)
\}.$
\begin{ex}\label{eq:adkxisquare}
${\rm ad}^{*k}_{\vert\xi\vert^2}(\tau)=( L_x+\Delta_x)^k\tau $ is a symbol of
order ${\rm ord}(\tau)+k$. 
\end{ex}
Here is another example of interest to us. 
\begin{ex}
\begin{eqnarray}\label{eq:logstar}
&&\{\log \vert \xi\vert^2, \tau\}_\star= \sum_{\vert \alpha  \vert =0 }^\infty\frac{(-i)^{\vert
  \alpha\vert}}{\alpha!}\partial_{\xi}^{\alpha}\log \vert \xi\vert^2
\partial_x^{\alpha}\tau (x,\xi)\\
&=&-2 i\, \sum_{j=1}^n \frac{\xi_j}{\vert \xi\vert^2}\partial_{x_ j}\tau
(x,\xi)-\sum_{i=1}^n \frac{1}{\vert \xi\vert^2}\partial_{x_i}^2\tau
(x,\xi)\nonumber+ 2 \, \sum_{i, j=1}^n 
\frac{\xi_i,\xi_j}{\vert\xi\vert^{4}}\, \partial_{x_ix_j}^2\tau
(x,\xi)
+\cdots\nonumber
\end{eqnarray}
\end{ex}
We now specialise to the algebra ${\cal F}S_{\rm cl}(U)$ of  polyhomogeneous formal symbols.  The
resolvent of a polyhomogeneous formal symbol $\sigma$ of order $a$
\begin{equation}\label{eq:resolvent} r_\star(\sigma, \lambda)= (\lambda- \sigma)^{\star^{-1}},
\end{equation}
solution of 
$ (\lambda- \sigma)\star r= 1$  has homogeneous components
$\sigma_{a-j}(r_\star(\sigma, \lambda))$ of degree $a-j$ in $(\xi,
\lambda^{\frac{1}{a}})$ defined inductively on $j$ by 
\begin{eqnarray}\label{eq:resolventstar}\sigma_{-a}(r_\star(\sigma, \lambda))&=&
  (\sigma_a-\lambda)^{-1}, \\ \sigma_{-a-j}(r_\star(\sigma,
  \lambda))
&=&-\sigma_{-a}(
r_\star(\sigma, \lambda))\sum_{k+l+\vert
  \alpha\vert=j,l<j}\frac{(- i)^{\vert \alpha\vert}}{\alpha!} \, D_\xi^\alpha
\sigma_{a-k}(\sigma) D_x^\alpha \sigma_{-a-l}(r_\star(\sigma, \lambda)).\nonumber
\end{eqnarray} 
\begin{defn} We call a  formal symbol $\sigma$ in ${\cal F}S_{\rm cl}(U)$
  admissible with spectral cut 
  $\theta$ if for every $(x, \xi)\in T^*U-\{0\}$  the
  leading symbol matrix  $\sigma^L(x, \xi)$ has no eigenvalue in a conical
  neighborhood of  the ray  $L_\theta=\{r
  e^{i\theta}, r\geq 0\}$. In particular, the symbol is elliptic.
\end{defn} 
The logarithm of an admissible formal polyhomogeneous symbol $\sigma$ is
defined by (see e.g. \cite{O1}): 
$$\log_\star(\sigma):=  \frac{i}{2\pi}\left(\partial_z
\int_\Gamma \lambda^{z}\,  (\lambda- \sigma)^{\star^{-1}}\,
d\lambda\right)_{\vert_{ z=0}},$$ for a contour $\Gamma$, which encloses the
eigenvalues of the leading symbol of $\sigma$.
The Campbell-Hausdorff formula for admissible formal  polyhomogeneous  symbols $\sigma$
and $\tau$ with scalar
leading symbols reads (see
Lemma 2.7 in \cite{O1}):
\begin{equation}\label{eq:CHsymb}
\log_\star(\sigma\star \tau )\sim \log_\star \sigma+\log_\star  \tau
+\sum_{k=2}^{\infty } C_\star^{(k)}(\log_\star\sigma , \log_\star\tau
  ),
\end{equation} where 
 $C_\star^{(k)}( \log_\star \sigma, \log_\star \tau)$ are  Lie monomials
 defined as in (\ref{eq:Ck}) replacing ${\rm ad}_P(Q)$ by 
 $$ {\rm ad}_p^\star (q):= \{p, q\}_\star:= p\star
  q-q\star p.$$
  The beginning of the expansion in equation (\ref{eq:CHop}) reads: 
\begin{eqnarray}\label{eq:CHexp}
\log_\star(\sigma\star \tau )&\sim & \log_\star \sigma+\log_\star  \tau
+\frac{1}{2} \{\log_\star \sigma, \log_\star \tau\}_\star \nonumber \\
&+&\frac{1}{12} \{\log_\star
\sigma, \{\log_\star \sigma, \log_\star \tau\}_\star\}_\star- \frac{1}{12} \{\log_\star
\tau, \{\log_\star \sigma, \log_\star \tau\}_\star\}_\star\nonumber\\
&-& \frac{1}{24} \{\log_\star \tau, \{\log_\star
\sigma, \{\log_\star \sigma, \log_\star \tau\}_\star\}_\star\}_\star\cdots
\end{eqnarray} 
\begin{rk}\label{rk:classicalbrackets} If $\tau$ is classical, then $C_\star^{(k)} (\log\vert  \xi\vert^2,
  \tau)$ is classical since the bracket $\{\log\vert \xi\vert^2, \sigma\}_\star$ with a
  classical symbol $\sigma$ is classical.
\end{rk}
\begin{rk}\label{rk:decreasing}  If $\tau$ has negative order, then
the order $\alpha_k$ of $C_\star^{(k)} (\log\vert  \xi\vert^2, \tau)$ is
negative and decreases with $k$. Indeed,    $\alpha_{k+1} $   corresponds
either to  the order of  $\left\{\log\vert
\xi\vert^2,  C_\star^{(k)} (\log\vert  \xi\vert^2, \tau)\right\}_\star$, which by
(\ref{eq:logstar}) is
$\alpha_k-1$ or to the order of  $\left\{\tau,  C_\star^{(k)} (\log\vert
  \xi\vert^2, \tau)\right\}_\star$ which is ${\rm ord}(\tau)+ \alpha_k$ and hence
smaller than $\alpha_k$. 
\end{rk}
\begin{thm} \label{thm:CHlogresdensity} The logarithmic  residue density (\ref{eq:logresdensity}) of        a
  generalised Laplacian  $Q$ on $E$
 is a finite sum of residue densities  of classical symbols: 
\begin{eqnarray}\label{eq:CHlogresdensity} {\rm res}_x(\log Q)&=& {\rm res}_x(\log_\star
(\vert \xi\vert^{-2}\star \sigma(Q)(x, \xi))\\
 &+&\sum_{j=1}^n \frac{(-1)^j}{(j+1)!} \, \sum_{k=2}^n{\rm res}_x\left( 
C_\star^{(k)} \left(\log \vert
  \xi\vert^2,\left(\vert \xi\vert^{-2}\star
  \sigma_{<2}(Q)(x, \xi)\right)^{* (j+1)}  \right)\right),\nonumber
  \\&=&\sum_{j=1}^n 
\frac{(-1)^j}{(j+1)!} \, {\rm res}_x\left(\left(\vert \xi\vert^{-2}\star
  \sigma_{<2}(Q)(x, \xi)\right)^{* (j+1)} \right)\nonumber\\
&+&\sum_{j=1}^n \frac{(-1)^j}{(j+1)!} \, \sum_{k=2}^n{\rm res}_x\left( 
C_\star^{(k)} \left(\log \vert
  \xi\vert^2,\left(\vert \xi\vert^{-2}\star
  \sigma_{<2}(Q)(x, \xi)\right)^{* (j+1)}  \right)\right),\nonumber
\end{eqnarray} where we have set  $\sigma(Q)(x, \xi)= \vert \xi\vert^2+ \sigma_{<2}( Q)(x,\xi)$.
\end{thm}
{\bf Proof:} We write $\sigma(Q)(x, \xi)=   \vert
\xi\vert^2\star \left(1+ \vert \xi\vert^{-2}\star \sigma_{<2}(Q)(x, \xi)
\right).$ Applying  the Campbell-Hausdorff formula  (\ref{eq:CHsymb})   to $\sigma= \vert \xi\vert^2$ and
 $\tau= 1+ \vert \xi\vert^{-2}\star
  \sigma_{<2}( Q)$ yields:
  \begin{eqnarray*}
\sigma(\log Q)(x, \xi)
&\sim& \log_\star  \sigma(Q)(x, \xi)\\
&\sim& 2\log\vert \xi\vert +\log_\star
(\vert \xi\vert^{-2}\star \sigma(Q)(x, \xi)\\
& +&\sum_{k=2}^\infty C_\star^{(k)} \left(\log \vert
  \xi\vert^2,\log_\star\left( 1+ \vert \xi\vert^{-2}\star
  \sigma_{<2}( Q)\right)\right)\\
&\sim& 2\log\vert \xi\vert +\sum_{j=1}^n \frac{(-1)^j}{(j+1)!} \, \left(\vert \xi\vert^{-2}\star
  \sigma_{<2}(Q)(x, \xi)\right)^{* (j+1)}\\
& +&\sum_{k=2}^\infty C_\star^{(k)} \left(\log \vert
  \xi\vert^2,\log_\star\left( 1+ \vert \xi\vert^{-2}\star
  \sigma_{<2}( Q)\right)\right).
\end{eqnarray*}
 This shows that $\log_\star \sigma(Q)-\log\vert \xi\vert^2$ is
  a clasical symbol as a consequence of  the fact that  the logarithm
$\log_\star  \tau 
\sim\sum_{j=1}^\infty \frac{(-1)^j}{(j+1)!} \, \left(\vert \xi\vert^{-2}\star
  \sigma_{<2}(Q)(x, \xi)\right)^{* (j+1)},$  is classical and hence  that  the
  corresponding Lie monomials are also classical  by Remark 
  \ref{rk:classicalbrackets}.   
Applying Remark \ref{rk:decreasing} to $\tau= \log_\star\left( 1+ \vert \xi\vert^{-2}\star
  \sigma_{<2}( Q)\right)$, which has negative order, shows that $ C_\star^{(k)} \left(\log \vert
  \xi\vert^2,\log_\star\left( 1+ \vert \xi\vert^{-2}\star
  \sigma_{<2}( Q)\right)\right)$ has order smaller than $ -k$. 
Since  the residue  vanishes on symbols of order smaller than $-n$ and $\left(\vert \xi\vert^{-2}\star
  \sigma_{<2}(Q)(x, \xi)\right)^{* (j+1)}$ has order no larger than $ -(j+1)$, implementing the
residue yields:
\begin{eqnarray*}
 {\rm res}_x(\log Q) 
&=&  {\rm res}_x\left(\log_\star\left(  
\vert \xi\vert^{-2}\star
  \sigma(Q)(x, \xi)\right) \right) \\
&+&\sum_{k=2}^{n}{\rm res}_x\left( C_\star^{(k)} \left(\log \vert
  \xi\vert^2\star\log_\star\left( 1+ \vert \xi\vert^{-2}\star
  \sigma_{<2}( Q)\right)\right)\right).
\end{eqnarray*}   
Replacing $\log_\star\left(  
\vert \xi\vert^{-2}\star
  \sigma(Q)(x, \xi)\right)= \log_\star\left( 1+ \vert \xi\vert^{-2}\star
  \sigma_{<2}( Q)\right)$ by its expansion yields the result.
\endsquare

\section{The logarithmic residue density via  a noncommutative Taylor expansion }
A noncommutative Taylor type formula provides an alternative way to express
  logarithmic
residue densities. We extend formulae for noncommutative Taylor expansions derived in \cite{P}
to formal polyhomogeneous symbols.   \\ Given an analytic function  $\phi(z)= \phi_0+\phi_1z+
\phi_2z^2+\cdots$ and  an admissible symbol $\sigma$ in ${\cal F}S_{\rm cl}(U)$
we write 
\begin{equation}\label{eq:Phi}
\Phi_\star(\sigma)= \frac{1}{2i\pi} \int_\Gamma r_\star(\lambda,\sigma)\, \phi(\lambda)\,
d\lambda.
\end{equation}
where the resolvent $r_\star(\lambda,\sigma)$ is defined by
(\ref{eq:resolvent}) and where $\Gamma$ is a contour which encloses the
eigenvalues of the leading symbol of $\sigma$. 
Applying this to the higher derivative $\phi^{(k)}$ yields:
\begin{equation}\label{eq:Phik}
\Phi_\star^{(k)}(\sigma)= \frac{1}{2i\pi} \int_\Gamma (\lambda- r_\star(\lambda,\sigma)\, \phi^{(k)}(\lambda)\,
d\lambda=\frac{k!}{2i\pi}\int_\Gamma 
\,(\lambda-\sigma)^{\star(-k-1)}
\, \phi(\lambda)\, d\lambda. 
\end{equation}
 If $\sigma= \vert \xi\vert^q+\sigma_{<q}$ with $\sigma_{<q}$ of order smaller
 than $q$,
 then the $\star$-resolvent reads:
 \begin{equation}\label{eq:resolvxa} r_\star (\lambda, \vert \xi\vert^q+\sigma_{<q})=r_\star(\lambda, \vert \xi\vert^q)+
   \sum_{n=1}^\infty  r_{\star n}(\lambda, \vert \xi\vert^q)
   \,(\sigma_{<q})^{\otimes n},
\end{equation}
where for symbols $\tau_1, \cdots, \tau_n$ in ${\cal F}S_{\rm cl}(U) $  we have set
\begin{eqnarray}\label{eq:resoln}
 && r_{\star n}(\lambda, \vert \xi\vert^q) (\tau_1\otimes \cdots \otimes
 \tau_n)\\
&=&  \sum_{\vert k\vert= 0}^\infty\frac{(k_1+\cdots+ k_{n}+n-1)!}{k! \,(k_1+1)(k_1+k_2+1)\cdots (k_1+\cdots+ k_{n-1}+n-1)}\nonumber\\
&\cdot& 
{\rm ad}_{\vert \xi\vert^q}^{(k_1)}(\tau_1)\star 
{\rm ad}_{\vert \xi\vert^q}^{(k_2)}(\tau_2)\star \cdots \star 
{\rm ad}_{\vert \xi\vert^q}^{(k_n)}(\tau_n)\,
(\lambda-\vert \xi\vert^q)^{-\vert k\vert-n-1},\nonumber
\end{eqnarray}
and $\vert k\vert= k_1+\cdots +k_n$ and 
  $k!= k_1!\cdots k_n!$.\\
Second quantised functionals  are defined  on tensor products of symbols  in terms of Cauchy integrals in analogy to ordinary
functionals on symbols (see  (\ref{eq:Phi})), but  by means of    quantised
resolvents
$r_{\star n}$ instead of the ordinary resolvent $r_\star$. 
\\ To  an   analytic function
$\phi(z)$ and to an admissible symbol $\sigma$   we assign a map called {\it
  second quantisation of $\Phi(x)$} defined on $\left({\cal F}S_{\rm cl}(U)\right)^{\otimes n}$  by
\begin{eqnarray*}
 \Phi_{\star\, n} ( \sigma):\quad \quad \quad\left({\cal F}S_{\rm cl}(U)\right)^{\otimes
   n} &\to& {\cal F}S_{\rm cl}(U) \\
\tau_1\otimes  \cdots\otimes  \tau_n&\mapsto&\frac{1}{2i\pi}\int_\Gamma
r_{\star n}(\lambda,
\sigma )(\tau_1\otimes \cdots\otimes \tau_n)\, \phi(\lambda)\, d\lambda.\\
\end{eqnarray*}
 One easily derives  the following noncommutative Taylor type formula from  (\ref{eq:resoln}):
\begin{equation}\label{eq:NCTaylorbis}
 \Phi_{\star \,n} (\sigma)\left(\tau_1\otimes \cdots \otimes \tau_n\right)=\sum_{\vert
  k\vert =0}^\infty \frac{{\rm ad}_\sigma^{\star k_1}(\tau_1)\star \cdots\star {\rm
    ad}_\sigma^{\star k_n}(\tau_n)}{k!\, (k_1+1)(k_1+k_2+2)\cdots (k_1+\cdots +k_{n}+n)}
\,  {\Phi^{(\vert k\vert+n)} (\sigma)}.
\end{equation} 
Applying (\ref{eq:resolvxa})   to $\sigma= \vert \xi\vert^2+\sigma_{<2}$ where $\sigma_{<2}$
  has order smaller than $2$ we have:
\begin{equation}\label{eq:NCTaylorPhi} \Phi_\star(\sigma)= \Phi_\star(\vert \xi\vert^q)+\sum_{p=1}^\infty
  \Phi_{\star\,p}(\vert \xi\vert^q)(\sigma_{<2}^{\otimes p})\end{equation}
Applying this to $\phi=\log$ yields
\begin{eqnarray}\label{eq:NCTaylorlog}&& \log_\star(\sigma)-\log_\star(\vert\xi\vert^2)\\
&=&\sum_{p=1}^\infty \sum_{\vert
  k\vert =0}^\infty (-1)^{\vert k\vert +p -1}\, (\vert k\vert+p-1)!\, \frac{
{\rm ad}_{\vert \xi\vert^2}^{\star k_1}(\sigma_{<2})\star {\rm ad}_{\vert \xi\vert^2}^{\star
  k_2}(\sigma_{<2})\cdots \star {\rm ad}_{\vert \xi\vert^2}^{\star
  k_p}(\sigma_{<2})}{k!\,(k_1+1)(k_1+k_2+2)\cdots (k_1+\cdots +k_{p}+p)}
\,  \vert  \xi\vert^{-2\,(\vert k\vert+p)}.\nonumber\end{eqnarray}
   Implementing the noncommutative residue finally leads to the
following formula for the   logarithmic residue density. 
\begin{thm}\label{thm:NCTaylorloggsymbol} The logarithmic residue density  of  a
  generalised Laplacian $Q$ on $E$ 
 is a finite sum of residues of
classical symbols: 
   \begin{eqnarray}\label{eq:NCTaylorlogsymbol}&& {\rm
       res}_x\left(\log (Q)\right)\\ 
&=&\sum_{p=1}^{n} \sum_{\vert
  k\vert =0}^{n-p} (-1)^{\vert k\vert +p -1}\, (\vert k\vert+p-1)!\,
\frac{{\rm res}_x\left((L_x+\Delta_x)^{k_1}( \sigma_{<2}(Q))\cdots (L_x+\Delta_x)^{k_p}( \sigma_{<2}(Q))\,  \vert \xi\vert^{-2(\vert k\vert+p)}\right)}{k!\,(k_1+1)(k_1+k_2+2)\cdots (k_1+\cdots +k_{p}+p)}
\nonumber,
\end{eqnarray}
where $k!:= k_1!\cdots k_p!$ and $\vert k\vert= k_1+\cdots+k_p$. Here we have
set $\sigma_{<2}(Q)(x, \xi):= \sigma(Q)(x, \xi)-\vert \xi\vert^2$ and as before, 
$L_x:= -2i  \sum_{a=1}^n \xi_a \partial_{x_a}$ and $\Delta_x:= -\sum_{a=1}^n
\partial^2_{x_a}$.
\end{thm}
{\bf Proof:} By (\ref{eq:NCTaylorlog}) combined with (\ref{eq:adkxisquare})  we have:
\begin{eqnarray*}
&&\sigma(\log (Q))(x, \xi)\\
&\sim&\sum_{p=1}^{\infty} \sum_{\vert
  k\vert =0}^\infty (-1)^{\vert k\vert +p -1}\, (\vert k\vert+p-1)!\,
\frac{(L_x+\Delta_x)^{k_1}( \sigma_{<2}(Q))\cdots (L_x+\Delta_x)^{k_p}(
  \sigma_{<2}(Q) )\,  \vert \xi\vert^{-2(\vert k\vert+p)}}{k!\,(k_1+1)(k_1+k_2+2)\cdots (k_1+\cdots +k_{p}+p)}
\nonumber,
\end{eqnarray*}
which is a formal power series of symbols $\sigma_k$ of  decreasing order
$-(\vert k\vert +p)$.  Since the noncommutative residue vanishes on symbols of
order smaller than $-n$, we have $\vert k\vert +p\leq n$ which implies that only terms  $p\leq n$ and
$\vert k\vert \leq n-p$ survive after applying the residue.\endsquare

 \section{The index  as a logarithmic (super-) residue}
Let us recall results of \cite{PS} and \cite{Sc1} (see also \cite{Sc2}).
 Let   $Q\in \Cl(M, E)$ be an
admissible (and hence invertible, see Appendix) classical pseudodifferential
operator of positive order $q$.
\\ For any {\it differential} operator  $A\in \Cl(M, E)$,   
the noncommutative residue density 
\begin{equation}\label{eq:residuedensitylog}{\rm res}_x(A\, \log Q)\, dx:= -\frac{1}{(2\pi)^n}\, \left(\int_{\vert\xi\vert=1} {\rm
  tr}\left(\sigma_{-n}(A\, \log Q)(x,\xi)\right)\,
d_S\xi\right)\, dx,
\end{equation}
is a globally defined $n$-form on $M$ (see \cite{O2} for the case $A=I$, \cite{PS} for
  the general case), which integrates over $M$ to the noncommutative residue:
\begin{equation}\label{eq:residuelog}{\rm res}(A\, \log Q):= -\frac{1}{(2\pi)^n}\,\int_M \left(\int_{\vert\xi\vert=1} {\rm
  tr}\left(\sigma_{-n}(A\, \log Q)(x,\xi)\right)\,
d_S\xi\right)\, dx.
\end{equation}
It furthermore relates to the  $Q$-weighted trace ${\rm Tr}^Q(A)$  of $A$ by (see \cite{Sc1} when $A=I$ and \cite{PS} for
  the general case)
 \begin{equation}\label{eq:trQreslog}
{\rm Tr}^Q(A):={\rm fp}_{z=0}{\rm TR}\left(A\,
  Q^{-z}\right)\nonumber
= -\frac{1}{q}\, {\rm
  res}(A\, \log_\theta Q),
\end{equation}
where ${\rm fp}_{z=0}$ stands for the finite part at $z=0$.
Here,  $d_S\xi$ is the  volume
 form on the unit sphere induced by the canonical
measure on $\R^n$, where $\sigma_{-n}$ stands for the positively homogeneous component of
degree $-n$ of a logpolyhomogeneous symbol $\sigma$.
\begin{rk}\label{rk:logsmoothing} Once checks that ${\rm res}(A\,\log (Q+R))=
  {\rm res}(A \,\log Q)$
for any smoothing operator $R$. 
\end{rk}
\begin{ex} Setting $A=I$ in the above
  corollary yields 
\begin{equation}\label{eq:zetaQ}
\zeta_Q(0)
= -\frac{1}{q}\, {\rm
  res}( \log Q),
\end{equation}
where  $\zeta_Q(z)$ is the zeta function associated to $Q$.
This corresponds to the logarithm $ {\rm
  res}( \log Q)=\log {\rm det}_{\rm res}(Q)$ of Scott's residue determinant
\cite{Sc1}.
\end{ex}
Let $E=E^+\oplus E^-$ be any $\Z_2$ graded vector bundle over $M$ and let  $D^+:\Cl(M,E_+)\to \Cl(M,E_-)$ be an elliptic operator in 
$\Cl\left(M,E_+^*\otimes E_-\right)$.
Its (formal) adjoint 
$D^-:=\left(D^+\right)^*: \Cl(M,E^-)\to \Cl(M,E^+)$ is an elliptic operator in
$\Cl\left(M,\left(E^-\right)^* \otimes
E_+\right)$ and $\Delta= \Delta^+\oplus \Delta^-$ with
$\Delta^+:= D^-\, D^+$, $\Delta^-:= D^+\, D^-$ are non-negative (formally)
self-adjoint elliptic operators. \\
 The following theorem which combines formulae due to   McKean and  Singer
  \cite{MS} and Seeley \cite{Se}, expresses the index of $D^+$:
$$ {\rm ind}(D^+):= {\rm dim}\left({\rm Ker}(D^+)\right)-  {\rm dim}\left({\rm
    Ker}(D^-)\right)$$  in terms of the superweighted trace of the
 identity. Let $\pi_\Delta$ denote the orthogonal projection onto the kernel
 of $\Delta$, which is finite dimensional as $M$ is compact.
\begin{thm}\label{thm:superresidue} The superresidue
$${\rm sres}\left( \log  (\Delta)\right):=  -\frac{1}{(2\pi)^n}\, \left(\int_{\vert\xi\vert=1} {\rm
  s tr}\left(\sigma_{-n,0}(\log  (\Delta))(x,\xi)\right)\,
d_S\xi\right)\, dx,$$
 is a globally defined $n$-form and   we have
\begin{equation}\label{eq:index}{\rm
  ind}(D^+)={\rm sTr}^{\Delta+\pi_\Delta}(I)=-\frac{1}{2\,{\rm
ord}(D)} {\rm sres}\left( \log  (\Delta+\pi_\Delta)\right),\end{equation}
where $\pi_\Delta$ is the orthogonal
projection onto the kernel of $\Delta$ and ${\rm
ord}(D)$ is the
order of $D$. Here str stand for the super trace on the graded fibres of $E$.
\end{thm}
\begin{rk} In view of Remark \ref{rk:logsmoothing}, one can drop the explicit
  mention of $\pi_\Delta$ and write  ${\rm sres}(\log \Delta)$ since  the
  projection $\pi_\Delta$ is  smoothing and the residue is invariant
  under translation by a smoothing operator.
\end{rk}
{\bf Proof:}  We first observe a property of the spectrum of $\Delta$:
$${\rm Spec}(\Delta^+)-\{0\}= {\rm Spec}(\Delta^-)-\{0\}.$$
Indeed, $$\Delta^+u_+=\lambda^+u_+\Rightarrow \Delta_- (D^+ u_+)= \lambda^+\,
D^+ u_+\quad \forall u_+\in \Ci(M, E_+)$$ so that an  eigenvalue $\lambda^+$
of $\Delta^+$ with eigenvector $u_+$ is  an
 eigenvalue of $\Delta^-$ with eigenvector $D^+u_+$ provided the latter does
 not vanish.  The converse holds similarly. 
\\ Let us denote by $\{\l_n^+, n\in \N\}$ the discrete set of eigenvalues of
$\Delta^+$  and by $\{\l_n^-, n\in \N\}$ the discrete set of eigenvalues of
$\Delta^-$.  For any complex number $z$:
\begin{eqnarray*}
{\rm sTr}\left((\Delta+\pi_\Delta)^{-z}\right)&=&\sum_{n\in\N} \left(\l_n^++ \delta_{\l_n^+}\right)^{-z}-
\sum_{n\in\N} \left(\l_n^-+ \delta_{\l_n^-}\right)^{-z}\\
&=&\sum_{\l_n^+\neq 0} \left(\l_n^+\right)^{-z}-
\sum_{\l_n^-\neq 0} \left(\l_n^-\right)^{-z} +{\rm dim} {\rm Ker}\Delta^+- {\rm dim} {\rm Ker}\Delta^-\\ 
&=& {\rm
  ind}(D^+).
\end{eqnarray*}
Taking the finite part at $z=0$ therefore yields:
$$
{\rm ind}(D^+)
= {\rm sTr}^{\Delta+\pi_\Delta}(I)= -\frac{1}{ 2} \
{\rm sres}
\left( \log (\Delta)\right).$$
\endsquare
\begin{ex} With the notations introduced at the beginning of the paper, for a
 Dirac operator $D_W^+: \Ci(M,
  S^+\otimes W )\to \Ci(M,  S^-\otimes W)$ on the $\Z_ 2$-graded spinor bundle $S=S^+\oplus S^-$
  over an even dimension spin manifold $M$ we have
\begin{equation}\label{eq:indexD}
{\rm ind }(D_W^+) 
= -\frac{1}{2}\, {\rm
  sres}\left( \log \left(D_W^2\right)\right)= -\frac{1}{2}\,
\int_M {\rm
  sres}_x\left( \log \left(D_W^2\right)\right)\, dx .
\end{equation}
\end{ex}
The remaining part of the paper deals with the computation of the logarithmic
density of the square $D^2$ of the Dirac operator $D$ acting on spinors. 
\section{A formula for the index via the Lichnerowicz formula}
   We first recall the  Lichnerowicz formula (see e.g. Theorem 3.52 of \cite{BGV})
    or equivalently the
general Bochner identity (see   Theorem 8.2
of \cite{LM}), which  relates the
square $D_W^2$ of the twisted Dirac operator $D_W$ 
 with the  Laplace-Beltrami operator
 \begin{equation}\label{eq:generalizedlaplacian}
\Delta^E=-{\rm
  tr}\left( \nabla^{T^*M\otimes E} \nabla^E\right)=-\sum_{i=1}^n
\left(\nabla^{T^*M\otimes E
}\nabla^E\right)_{e_i,e_i}=-\sum_{i=1}^n
 \left(\nabla^E_{e_i} \nabla^E_{e_i}-\nabla^E_{\nabla^{
     E}_{e_i} e_i}\right)
\end{equation}
 associated with
the superconnection $\nabla^E$ on $E$, where $\nabla^{T^*M\otimes E}$ is the connection induced on the tensor product bundle $T^*M\otimes E$ by
the Levi-Civita connection on $M$ and the connection $\nabla^E$ on $E$.
Here $\{e_i, i=1, \cdots, n\}$ is a local orthonormal tangent frame.
\begin{prop}
\begin{equation}\label{eq:Bochner}D_W^2=   \Delta^E+ R^E\nonumber\\
 =   \Delta^E+ R^W+\frac{r_M}{4},
  \end{equation}
where $r_M$ stands for the scalar curvature on $M$ and 
\begin{equation}\label{eq:RE}
R^E:=  \sum_{i<j} c(e_i) \, c(e_j)\,  \left( \nabla^{E}\right)^2_{e_i, e_j};\quad  R^W:=\sum_{i<j}  c(e_i)\,
 c(e_j)\,
 \left( \nabla^{W}\right)^2_{e_i, e_j}.\end{equation}
In particular, for a flat  auxillary bundle we have: 
$$D_W^2= \Delta_M +\frac{r_M}{4},$$
where $\Delta_M$ is the Laplace-Beltrami operator on the Riemannian manifold
$M$.
 
\end{prop}
{\bf Proof:} We choose a local orthonormal tangent frame $\{e_i, i=1, \cdots, n\}$
at point $x\in M$ such that $\left(\nabla^E_{e_i}\right)_x=0$ for all
$i\in \{1, \cdots, n\}$. Since $D_W= \sum_{i=1}^n c(e_i) \, \nabla^{E
}_{e_i}$, at  that point $x$ we have:
\begin{eqnarray*}
D_W^2&=&  \sum_{i,j=1}^n c(e_i) \, \nabla^{E
}_{e_i}\,c(e_j) \, \nabla^{E
}_{e_j}\\
&=& \sum_{i,j=1}^n c(e_i) \, c(e_j) \,\left[\left(\nabla^{E
}\right)^2_{e_i,e_j}+  \nabla^{E
}_{\nabla_{e_i} e_j}\right]\\
&=& -\sum_{i=1}^n \left(\nabla^{E
}\right)^2_{e_i,e_i}+\sum_{i<j} c(e_i) \, c(e_j) \,\left[\left(\nabla^{E
}\right)^2_{e_i,e_j}- \left(\nabla^{E  
}\right)^2_{e_j,e_i}\right]\\
&=& \Delta^E+\sum_{i<j} c(e_i) \, c(e_j)\, \left(\nabla^{E
}\right)^2_{e_i, e_j}\\
&=& \Delta^E+R^E.
\end{eqnarray*}
The curvature term $\left(\nabla^E\right)^2\in \Omega^2(M,{\rm End} (E))$
decomposes as $ \left(\nabla^E\right)^2= \left( \nabla^S\right)^2
\otimes 1+ 1\otimes  \left( \nabla^{W}\right)^2$ so that
$R^E= \sum_{i<j}
c(e_i) \, c(e_j)\, \left( \nabla^S\right)^2_{e_i, e_j}+ R^W.$  A careful computation (see e.g. the proof of
Theorem 3.52 in \cite{BGV}) shows that $ \sum_{i<j}
c(e_i) \, c(e_j)\, \left( \nabla^S\right)^2_{e_i, e_j}= \frac{r_M}{4}.$
 \endsquare\\
 Combining the Lichnerowicz formula with the Campbell-Hausdorff formula yields
 a formula for the index. 
\begin{thm}\label{thm:indtwisted} In even dimension $n=2p$, 
\begin{eqnarray}\label{eq:indres}
{\rm  ind}(D_W^+)&=&-\frac{1}{2} {\rm sres}(\log (D_W^2))\\
&=&-\frac{1}{2} {\rm sres}(\log
  (\Delta^E ))+\sum_{k=1}^{n-1} \frac{(-1)^{k}}{2\, k}{\rm
    sres}\left( \left[ (\Delta^E)^{-1}\, R^E\right]^k
  \right).\nonumber
\end{eqnarray}
Inside the residue we write for short $ (\Delta^E)^{-1}$ instead of $ (\Delta^E+\pi_\Delta)^{-1}$
since the residue is insensitive to the smoothing operator $\pi_\Delta$.  
\end{thm} 
{\bf Proof:}  By equation (\ref{eq:Bochner})
\begin{eqnarray*} D_W^2+\pi_{D_W^2}&=&  \Delta^E+\pi_{\Delta_E}+ R^E+
  \pi_{D_W^2}-\pi_{\Delta^E}\\
&=& \left(\Delta^E+\pi_{\Delta^E}\right)\, \left(1 + \left(\Delta^E+\pi_{\Delta^E}\right)^{-1}\left(R^E+
  \pi_{D^2}-\pi_{\Delta^E}\right)\right),
\end{eqnarray*}
so that by (\ref{eq:resCH}), we get:
\begin{eqnarray*}&{}&{\rm sres}\left(\log( D_W^2)\right)\\
&=& 
{\rm sres}\left(\log \left(\Delta^E \right)\right)+ {\rm sres}\left(\log \left(1 + \left(\Delta^E \right)^{-1}\left(R^E
  \right)\right)\right)\\
&=& {\rm sres}\left(\log
  \left(\Delta^E \right)\right)+\sum_{k=1}^\infty
\frac{(-1)^{k+1}}{k}\,  {\rm sres}\left(\left[\left(\Delta^E \right)^{-1}\left(R^E
  \right)\right]^k\right)\\
&=& {\rm sres}\left(\log
  \left(\Delta^E \right)\right)+\sum_{k=1}^\infty
\frac{(-1)^{k+1}}{k}\,  {\rm
  sres}\left(\left[\left(\Delta^E \right)^{-1} \,R^E\right]^k\right).
\end{eqnarray*}
Here we  used the fact that the noncommutative residue
vanishes on smoothing operators. Also,  for an operator $B\in \Cl(M, E)$ with
negative order, we have ${\rm sres}\left(\log(1+B)\right)= \sum_{k=1}^\infty
\frac{(-1)^{k+1}}{k}{\rm sres}(B^k),$ which  is actually a finite sum since the residue
vanishes for operators of order smaller than minus the dimension of the
underlying manifold. Since $B= \left(\Delta^E+\pi_\Delta \right)^{-1}\,
R^E$ has order $-2$, the sum stops at $p= n/2$.
 \endsquare

\section{The Atiyah-Singer index theorem for a twisted Dirac operator on a flat
space}
We derive  the Atiyah-Singer  index formula for  a twisted Dirac operator on a flat space  from 
(\ref{eq:indres}). We use the
 notations introduced at the beginning of the paper.
 Denoting by $\partial_i + A_i$ the components of the connection $\nabla^E$ 
 in a given local trivialization of $E$ and local coordinates $x_i$ on $M$
 such that
 the metric Christoffel symbols vanish,  according to the Lichnerowicz formula
 we have
 $${D_W}^2 = \Delta^E + R^E= \sum_i (\partial_i + A_i)^2  + 
 \sum_{i<j} \gamma_i\gamma_j  F_{ij},$$ where we have set
 $F_{ij} = \partial_i A_j -\partial_j A_i + [A_i, A_j]$ to be
 the 2-form components of the curvature $\left(\nabla^E\right)^2.$ 
 \begin{thm}\label{thm:flattorus} If  the Riemann metric on $M$ is flat, then
$${\rm ind}(D_W^+)=  \int_M {\rm
  \tr}\left(e^{i\, \frac{F}{2\pi}}\right).$$
\end{thm}
{\bf Proof:}  As before, $n=2p$ stands for the dimension of $M$.
By Theorem \ref{thm:indtwisted} we have
$$ {\rm sres}\log({D_W}^2) = {\rm sres}(\log \Delta^E) +\sum_{k=1}^{n-1} \frac{(-1)^{k+1}}{ k}{\rm
    sres}\left( \left[ (\Delta^E)^{-1}\, R^E\right]^k
  \right).$$
The first term on the r.h.s. vanishes. Indeed, at a given point $x\in M$ and for
fixed $\xi\in T_x^*M$,  the $-n$-th homogeneous component of the symbol
 $\sigma(\log \Delta^E)(x, \xi)$ of
 $\log \Delta^E$ is an   endomorphism
 of the fibre $W_x$ of the auxillary vector bundle $W$.  By  
 (\ref{eq:strgamma})  the fibrewise
supertrace therefore  vanishes on  the $-n$-th homogeneous component
 of the symbol and hence   so does the residue density 
 ${\rm sres}_x(\log \Delta^E)\, dx$. Thus ${\rm
 sres}(\log \Delta^E)=0$.
We now investigate the second term on the r.h.s. On the one hand,  all the 
expressions ${\rm sres} 
\log\left(\left[ (\Delta^E)^{-1}\, R^E\right]^k \right)$  inside the sum 
vanish for   $k > p$, for
  the operators
 $\left[ (\Delta^E)^{-1}\, R^E\right]^k$  being  of order smaller than $-n$,
 their residues vanish.\\ 
   On the other hand, the expressions inside the sum also vanish for $k<p$ . 
   Indeed, 
  at a point $x\in M$ and for fixed $\xi\in T_x^*M$, the symbols in the
  variables $(x, \xi)$ inside the
   residues are of the form  $M \,\gamma_1\, \gamma_2\cdots\gamma_{k}$ 
    for some matrix $M\in {\rm End}(W_x)$
 and   sets
 $\{i_1, \cdots, i_{k}\}$  strictly smaller than $ \{1, \cdots, n\}$.  Their
 supertrace whcih arise inside the superresidue, therefore vanishes by (\ref{eq:strgamma}).\\
 The remaining      $k=p$ term in the sum, which corresponds to the residue of an operator
 of order $-n$, only involves the leading symbol 
 $\sigma_L(\Delta^E)=\vert \xi\vert^2$ 
of $\Delta^E$. Thus we obtain
\begin{eqnarray*}
{\sres } \log({D_W}^2) &=& -\frac{(-1)^{p}}{ p}{\rm
    sres}\left( \left[ (\Delta^E +\pi_{\Delta^E})^{-1}\, R^E\right]^p \right)\\
     &=& -\frac{(-1)^{p} }{p} {\sres}\left(\vert \xi\vert^{-n} 
     \left({\rm tr}(R^E) \right)^{p}\right)\\ 
    &=& -\frac{(-1)^{p} }{p} \int_{M} 
    {\sres}_x\left(\vert \xi\vert^{-n} \,\left(\sum_{i< j}\gamma_i\, \gamma_j 
    F_{ij}\right)^p\right) \, dx\quad {\rm by}\quad (\ref{eq:RE}) \\
     &=& -\frac{(-1)^{p}\, 2^{p} \, {\rm Vol}(S^{n-1})}{(2\pi)^n\, p} \int_{M} 
   {\rm str} \left( \sum_{i,j}\sigma_{ij}\,  F_{ij}\right)^p \, dx\quad 
   {\rm by}\quad (\ref{eq:sigmaij})\quad
     {\rm and}\quad
     (\ref{eq:resintro})  \\ 
     &=& -2 \frac{i^p }{(4\, \pi)^p\, p!}
    \sum_{\tau \in \Sigma_n} (-1)^{\vert \tau\vert}  \int_{M}   
    F_{\tau(1)\tau(2)}\cdots  F_{\tau(n-1)\tau(n)}  \, dx 
     \quad {\rm by}\quad (\ref{eq:strprodsigma})\quad {\rm and}\quad 
     (\ref{eq:volsphere}) \\  
     &=& -2  \frac{i^p }{(2\, \pi)^p\, p!}\int_{M} 
    {\rm tr}\left( F^{\wedge p} \right) \, dx  \\
     &=& -2\int_M {\rm tr}\left(e^{i\, \frac{F}{2\pi}}\right).
\end{eqnarray*}Here, we have used
\begin{equation}\label{eq:volsphere}
{\rm vol}(S^{n-1})= \frac{2\pi^{\frac{n}{2}}}{\Gamma\left(\frac{n}{2}\right)} 
= \frac{2\pi^p}{ (p-1)!}.\end{equation}  
\endsquare

\section{The curvature tensor in normal coordinates}
  We recall a few properties of the curvature in a normal local
  coordinate system, i.e.,  a coordinate system defined by the exponential map
  at a point, so that rays emanating from the origin in the tangent space at a
  point are mapped to geodesics on the manifold emanating from this point. Let us   recall
that in Riemannian normal coordinates (see e.g. Proposition 1.28 in \cite{BGV}), \begin{equation}\label{eq:metricnormal}g_{ij}= \partial_{ij}-\frac{1}{3}R_{ikjl}x^l x^k+\sum_{\vert \alpha\vert
  \geq 3}\partial_\alpha g_{ij}\frac{x^\alpha}{\alpha!}\end{equation} 
 \begin{lem}\label{lem:Rsigma}  We have
$$\left(R_{iajk}+  R_{ikja}\right)\sigma_{kj}=\frac{3}{2}
R_{iajk}\sigma_{kj},$$
where $\sigma_{ij}$ was defined in (\ref{eq:sigmaij}).
\end{lem}
{\bf Proof:} 
Using the first Bianchi identity
$$R_{[ijk]l}=0$$ we  write $R_{ijka}= -R_{kija}-R_{jkia}$, which combined with the
antisymmetry of $\sigma_{ij}$ in $i$ and $j$ and the (anti)symmetry properties
of the curvature tensor $R_{ijkl}= -R_{jikl}=-R_{ijlk}= R_{klij}$ yields:
\begin{eqnarray*}
\left(R_{iajk}+  R_{ikja}\right)\sigma_{kj}
&=&R_{iajk}\sigma_{kj}+ R_{ijka}\sigma_{jk}\\
&=& \left(R_{iajk}+ R_{kija}+R_{jkia}        \right)\sigma_{kj}\\
&=& \left(2 R_{iajk}- R_{ikja}       \right)\sigma_{kj}.
\end{eqnarray*}
 Consequently, 
\begin{eqnarray*}
\left(R_{iajk}+  R_{ikja}\right)\sigma_{kj}&=&\frac{1}{2}\left[\left(R_{iajk}+
  R_{ikja}\right)\sigma_{kj}+ \left(2 R_{iajk}- R_{ikja}
\right)\sigma_{kj}\right]\\
&=& \frac{3}{2} R_{iajk}\sigma_{kj}
\end{eqnarray*}\endsquare
\begin{prop}\label{prop:partialGammasigma}
At the center of  a normal coordinate system, we have
\begin{equation}\label{eq:partialGammasigma}\partial_a
  \Gamma_{ij}^k\sigma_{kj}=\frac{1}{2} R_{jkia}\sigma_{kj}
\end{equation} so that 
$\partial_i \Gamma_{ij}^k\sigma_{kj}=0$. 
\end{prop}
{\bf Proof:} 
 By (\ref{eq:metricnormal}) the Christoffel symbols $\Gamma_{ij}^k= \frac{1}{2}g^{kl}\left( \partial_j g_{il} +
   \partial_i g_{jl}- \partial_l g_{ij}\right)$
vanish at the center of the normal coordinate symstem, where we have:
$$\partial_a \Gamma_{ij}^k=\frac{1}{3}\left(R_{iajk}+ R_{jaik}\right).$$
Indeed, differentiating  (\ref{eq:metricnormal}) twice yields
\begin{eqnarray*}
\partial_a\Gamma_{ij}^k&=& \frac{1}{2}\, \delta^{kl} \, \left( \partial_a
  \partial_j g_{il}+ \partial_a \partial_i g_{jl}-\partial_a \partial_l
  g_{ij}\right)\\
&=& -\frac{1}{6}\,\left(R_{iakj}+ R_{ijka}+ R_{jaki}+
  R_{jika}-R_{iajk}-R_{ikja}\right)\\
&=& \frac{1}{3}\, \left(R_{iajk}+ R_{jaik}\right).
\end{eqnarray*}
It follows from Lemma \ref{lem:Rsigma} that
$\partial_a \Gamma_{ij}^k\, \sigma_{kj}=\frac{1}{2}R_{iajk}\sigma_{kj}.$ and
hence in particular that 
$\partial_i \Gamma_{ij}^k\sigma_{kj}=\frac{1}{2}R_{iajk}\sigma_{kj}=0.$
 at the center of the normal
coordinate system.
\endsquare\\  The following result is useful to compute the index. 
\begin{prop}\label{prop:srespartialGamma} In four dimensions we have:
\begin{equation}\label{eq:srespartialGamma4} {\rm sres}_x\left( \vert \xi\vert^{-4} \,\partial_{x_{a}}\Gamma_{ij}^{k}\partial_{x_{a}} \Gamma_{im}^{n}\,
\sigma_{kj}\sigma_{nm}\right)\, dx
= \frac{1}{32\,\pi^{2}}\,{\rm tr}(R\wedge R),
\end{equation} and 
\begin{equation}\label{eq:srespartialGamma4bis} {\rm
    sres}_x\left(\frac{\xi_{a}\xi_{b}}{   \vert
      \xi\vert^{6}}
    \,\partial_{x_{a}}\Gamma_{ij}^{k}\partial_{x_{b}}\Gamma_{im}^{n}\right)\,
  dx
= \frac{1}{ 4\times 32\, \pi^{2}}\,{\rm tr}(R\wedge R).
\end{equation}
\end{prop}
{\bf Proof:} The result in four dimensionsis a consequence of the following   formula
in $n=2p=4q$ dimensions. At the center of   a normal coordinate
  system we show that:
\begin{eqnarray}\label{eq:prodsigmatau}
&&{\rm sres}_x\left( \vert \xi\vert^{-2p} \,\partial_{x_{a_1}}\Gamma_{i_1j_1}^{k_1}\partial_{x_{a_1}} \Gamma_{i_1m_1}^{n_1}\cdots\partial_{x_{a_p}}\Gamma_{i_qj_q}^{k_q}\partial_{x_{a_q}} \Gamma_{i_qm_q}^{n_q}\,
\sigma_{k_1j_1}\sigma_{n_1m_1}\cdots\sigma_{k_qj_q}\sigma_{n_qm_q}\right)\, dx
\nonumber\\
&=& \frac{ 1}{\Gamma(p)\, 2^{3p-1}\, \pi^p}\,\left({\rm tr}(R\wedge R)\right)^q.
\end{eqnarray}
The proof follows from combining   (\ref{eq:partialGammasigma})
with (\ref{eq:residuedensitylog}) and  the  formula for the volume of the unit
sphere $S^{n-1}$ in $n$  dimensions given by (\ref{eq:volsphere}):
\begin{eqnarray*}
&&{\rm sres}_x\left( \vert \xi\vert^{-2p} \,\partial_{x_{a_1}}\Gamma_{i_1j_1}^{k_1}\partial_{x_{a_1}} \Gamma_{i_1m_1}^{n_1}\cdots\partial_{x_{a_q}}\Gamma_{i_qj_q}^{k_q}\partial_{x_{a_q}} \Gamma_{i_qm_q}^{n_q}\,\sigma_{k_1j_1}\sigma_{n_1m_1}\cdots\sigma_{k_qj_q}\sigma_{n_qm_q}
\right)\\
&=&\frac{1}{2^{2q}}{\rm sres}_x\left( \vert \xi\vert^{-2p}
  \,R_{i_1a_1j_1k_1}\, R_{i_1a_1 m_1n_1}\cdots R_{i_qa_qj_qk_q}\, R_{i_qa_q m_qn_q}\,\sigma_{k_1j_1}\sigma_{n_1m_1}\cdots\sigma_{k_qj_q}\sigma_{n_qm_q} 
\right)\quad {\rm  by}\quad  (\ref{eq:partialGammasigma})\\
&=&\frac{(-i)^p}{4^p}\, \sum_{\tau \in \Sigma_{n}}(-1)^{\vert \tau\vert}\,  {\rm res}_x\left( \vert \xi\vert^{-2p}
  \,R_{i_1a_1\tau(1)\tau(2)}\,
 R_{i_1a_1\tau(3)\tau(4)}\,\cdots R_{i_qa_q\tau(n-3)\tau(n-2)}\, R_{i_qa_q
   \tau(n-1)\tau(n)}\right)\\
&\quad& {\rm  by}\quad (\ref{eq:strgamma})  \quad{\rm and}\quad 
(\ref{eq:strprodsigma})
\\
&=&\frac{ 1}{\Gamma(p)\, 2^{2n-1}\, \pi^p}\, 
\sum_{\tau \in \Sigma_{n}}(-1)^{\vert \tau\vert}\,  
  R_{i_1a_1\tau(1)\tau(2)}\,
 R_{a_1i_1\tau(3)\tau(4)}\,\cdots R_{i_qa_q\tau(n-3)\tau(n-2)}\, R_{a_qi_q
   \tau(n-1)\tau(n)}.
\end{eqnarray*}
  On the other hand we have
$$\frac{1}{2^p}\sum_{\tau \in \Sigma_{n}}(-1)^{\vert \tau\vert}\,  
  R_{i_1a_1\tau(1)\tau(2)}\,
 R_{a_1i_1\tau(3)\tau(4)}\,\cdots R_{i_qa_q\tau(n-3)\tau(n-2)}\, R_{a_qi_q
   \tau(n-1)\tau(n)}\, dx= \left({\rm tr}(R\wedge R)\right)^q,$$ which yields
\begin{eqnarray*}
&&{\rm sres}_x\left( \vert \xi\vert^{-2p} \,\partial_{x_{a_1}}\Gamma_{i_1j_1}^{k_1}\partial_{x_{a_1}} \Gamma_{i_1m_1}^{n_1}\cdots\partial_{x_{a_q}}\Gamma_{i_qj_q}^{k_q}\partial_{x_{a_q}} \Gamma_{i_qm_q}^{n_q}\,\sigma_{k_1j_1}\sigma_{n_1m_1}\cdots\sigma_{k_qj_q}\sigma_{n_qm_q}
\right)\, dx\\
&=&\frac{ 1}{\Gamma(p)\, 2^{3p-1}\, \pi^p}\,\left({\rm tr}(R\wedge R)\right)^q.
\end{eqnarray*} This proves the first part of the statemnt.
We prove the second part similarly. We first
observe that using the symmetries of the sphere we
have\begin{equation}\label{eq:intsphere} \int_{\vert \xi\vert=1} \frac{\xi_{i}\xi_j}{\vert
  \xi\vert^{n+2}}\, d\xi=  \delta_{i-j}\int_{\vert \xi\vert=1} \frac{\xi_{i}^2}{\vert
  \xi\vert^{n+2}}\, d\xi= \frac{1}{n} \delta_{i-j}\int_{\vert \xi\vert=1}
\frac{\sum_{i=1}^n \xi_i^2}{\vert
  \xi\vert^{n+2}}\, d\xi=\frac{\delta_{i-j}}{n} \,
\frac{2\pi^{\frac{n}{2}}}{\Gamma\left(\frac{n}{2}\right)}\end{equation}
so that  in four dimensions we get:
$${\rm sres}_x\left(\frac{\xi_{a}\xi_{b}}{   \vert \xi\vert^{6}}
  \,\partial_{x_{a}}\Gamma_{ij}^{k}\partial_{x_{b}} \Gamma_{im}^{n}\right)\,
dx= \frac{1}{4\times 32\, \pi^2} \,{\rm tr}\left(R\wedge R\right).$$\endsquare

\section{The Atiyah-Singer index formula in dimension $4$}
The square  of a Dirac operator $D$ acting on pure spinors is the prototype of
a generalised Laplacian. In local coordinates, its symbol reads:
$$\sigma(D^2)= \vert \xi\vert^2+ \Gamma_{ij}^k \sigma_{kj}\xi_i+
\partial_i \Gamma_{ij}^k \sigma_{kj} + \Gamma_{ij}^k \Gamma_{im}^n\,\sigma_{kj}
\sigma_{nm}+s= \vert \xi\vert^2+\sigma_{<2}(D^2),$$
where 
\begin{equation}\label{eq:sigmasmaller2} \sigma_{<2}(D^2):=\Gamma_{ij}^k \sigma_{kj}\xi_i+\partial_i \Gamma_{ij}^k \sigma_{kj} + \Gamma_{ij}^k\Gamma_{im}^n\, \sigma_{kj}
\sigma_{nm}+s, 
\end{equation}
and where $s$ stands for the scalar curvature.
\\ 
We use  (\ref{eq:indexD}) to compute the index of $D^+$: 
$${\rm ind}(D^+)= -\frac{1}{2}\,\int_M {\rm sres}_x(\log D^2)\, dx,$$ 
 in terms of the logarithmic (super) residue density,  which we explicitly derive in
four dimensions.\\
Since the residue does not depend on the choice of local coordinate, we choose
to derive
the residue density in a normal coordinate system.
We therefore need to compute:
$${\rm sres}_x(\log D^2)= {\rm sres}_x\left(\log_\star (\vert
\xi\vert^2+\sigma_{<2}(D^2))\right),$$ where $\log_\star$ is the logarithm on
symbols. There are at least two methods to compute the logarithm of $\vert
\xi\vert^2+\sigma_{<2}(D^2)$,  the 
 Campbell-Hausdorff 
formula (\ref{eq:CHsymb}) and  a Taylor type formula as in
 (\ref{eq:NCTaylorlogsymbol}).The first method, which in four dimensions
 reads \footnote{ As before, here  $\sigma_{< k}$ stands for the part of the symbol $\sigma$ of
order smaller than $k$.}:
 \begin{eqnarray}\label{eq:sreslogD2} &&{\rm sres}_x(\log(D^2))
 \nonumber\\
 &=&
 {\rm sres}_x
 \left(\log_\star (\vert \xi\vert^{-2}\star\sigma(D^2)) \right)\nonumber\\
&+&\frac{1}{2}\,\sum_{j=0}^2 \frac{(-1)^j}{(j+1)!} \,{\rm sres}_x\left( \left( 
\log \vert
  \xi\vert^2\star \left(  \vert
  \xi\vert^{-2}\star \sigma_{<2}(D^2)(x, \xi)\right)^{*(j+1)} \right)_{<-j-1}\right)\nonumber\\
&+&\frac{1}{12}\,\sum_{j=0}^1 \frac{(-1)^j}{(j+1)!} \,{\rm sres}_x
\left(\left(\log \vert
  \xi\vert^2\star  \left(  \log \vert
  \xi\vert^2\star\left( \vert
  \xi\vert^{-2}\star \sigma_{<2}(D^2)(x, \xi)\right)^{*(j+1)}\right)_{<-j-1}
  \right)_{<-j-2}\right)\nonumber\\
  &-& \frac{1}{12}\, {\rm sres}_x\left( 
  \{
\vert
  \xi\vert^{-2}\star \sigma_{<2}(D^2),\left(\log \vert
  \xi\vert^2\star \left(  \vert
  \xi\vert^{-2}\star \sigma_{<2}(D^2)(x, \xi) \right) \right)_{<-1}\}_\star
  \right),
 \end{eqnarray} and therfore requires computing the various terms in the
 above sums, 
is lengthier than the second method, which we adopt here. 
 \\
Replacing the residue by a super residue in  (\ref{eq:NCTaylorlogsymbol})  yields the following
description of the logarithmic superresidue density of $D^2$, in which the
sum over $p$ reduces to one term.  
\begin{prop} \label{prop:viaNCT} The logarithmic super residue density  of the squared Dirac operator
is a finite sum of super residues of classical symbols:
\begin{eqnarray}&& {\rm sres}_x(\log(D^2))\nonumber\\
&=& \sum_{ \vert k\vert = \frac{n}{4}, k_i\in\{1,2\}}^{\frac{n}{2}} \,
\frac{(-1)^{\vert k\vert +q-1}\, (\vert k\vert +q-1)!}{k!\,(k_1+1)(k_1+k_2+2)\cdots (k_1+\cdots +k_{p}+p)}\\
&\times &{\rm sres}_x\left((L_x+\Delta_x)^{k_1}(\sigma_{<2}(D^2))\cdots
 (L_x+\Delta_x)^{k_q}(\sigma_{<2}(D^2))\,  \vert \xi\vert^{-2(\vert k\vert+q)}\right)
\nonumber
\end{eqnarray}
  where as before,  we have set  $q=\frac{n}{4}$.
It is of the form:
$$ {\rm sres}_x(\log(D^2))= \sum_{s+t=q} \alpha_{s,t}\, {\rm
  sres}_x\left( (L_x^2\sigma_{<2}(D^2))^s\, (\Delta_x\sigma_{<2}(D^2))^t\,  \vert \xi\vert^{-2(3s+2t)}\right)$$
with 
$\Delta_x\sigma_{<2}(D^2)$ and $L_x^2\sigma_{<2}(D^2)$  contributing  respectively by
\begin{equation}\label{eq:Deltatau}\Delta_x\left(\Gamma_{ij}^k\Gamma_{lm}^n\right)
\sigma_{jk}\sigma_{mn}=-\frac{1}{2}
R_{jkia}R_{nmia}\,\sigma_{kj}\sigma_{nm}\end{equation} and
\begin{equation}\label{eq:Lsquare}L_x^2\left(\Gamma_{ij}^k\Gamma_{lm}^n\right)
\sigma_{jk}\sigma_{mn}=-  \, R_{jkia} R_{mnib}\,
\sigma_{jk} \sigma_{mn}\, \xi_a\,\xi_b.\end{equation}  
\end{prop} 
{\bf Proof:}
Applying (\ref{eq:NCTaylorlogsymbol}) to $Q=D^2$ yields an expression which
involves terms $(L_x+\Delta_x)^{k_i}(\sigma_{<2}(D^2))$, each of which differentiates
$\sigma_{<2}(D^2)$ at least $k_i$ times. We  have $k_i\leq 2$; indeed,  ${\rm
  sres}_x(\log_\star(D^2))$ being proportional to a Pontryagin form, it only
involves curvature terms so that 
only  first order derivatives of the Christoffel symbols can arise. In
view of the product term $\Gamma_{ij}^k \Gamma_{im}^n$, it can involve at most partial
differential operators of
order two.  \\  Since the superresidue density
${\rm sres}_x(\log D^2)\, dx$ is proportional to a Pontryagin form, there is
no contribution from the scalar curvature, so that  terms $\sigma_{<2}(D^2)$ corresponding
to zero powers $k_i$ do not contribute. \\
Let us analyse the contribution of terms involving powers $k_i=1$ i.e.,  expressions of the type
$(L_x+\Delta_x)\sigma_{<2}(D^2)$. In view of (\ref{eq:sigmasmaller2}) the terms $\Delta_x\sigma_{<2}(D^2)$ can only contribute by
\begin{eqnarray}\label{eq:Deltax}
\Delta_x \left(\Gamma_{ij}^k \Gamma_{im}^n\right)\sigma_{kj}\sigma_{nm}&=& -2\partial_a \Gamma_{ij}^k
\partial_a\Gamma_{im}^n \,\sigma_{kj}\sigma_{nm}\nonumber\\
&=& -\frac{1}{2} R_{jkia}R_{nmia}\,\sigma_{kj}\sigma_{nm}.
\end{eqnarray}  Let us now  see how the terms
$L_x\sigma_{<2}(D^2)$  contribute. 
We have $$L_x\sigma_{<2}(D^2)=-2i\left(\partial_a\Gamma_{ij}^k \sigma_{kj} \xi_i\xi_a+
  \partial_a\partial_i \Gamma_{ij}^k\sigma_{kj}\xi_a+(\partial_a \Gamma_{ij}^k
  \, \Gamma_{lm}^n+ \Gamma_{ij}^k\partial_a \Gamma_{im}^n)
  \sigma_{kj}\sigma_{nm}\xi_a+ \partial_a s\xi_a\right),$$
which at the center of a normal coordinate system reads:
$$L_x\sigma_{<2}(D^2)=-2i\left(\partial_a\Gamma_{ij}^k \sigma_{kj} \xi_i\xi_a+
  \partial_a\partial_i \Gamma_{ij}^k\sigma_{kj}\xi_a+\partial_a
  s\xi_a\right).$$
The only possible contribution  can come from \begin{equation}\label{eq:Lx}L_x
\left( \Gamma_{ij}^k\sigma_{kj}\xi_i\right)=-2i\,\partial_a\Gamma_{ij}^k
\sigma_{kj} \xi_i\xi_a = -i\ R_{jkia}\sigma_{kj} \,\xi_i\,\xi_a,
\end{equation}
which vanishes by antisymmetry of $R$. There is therefore no contribution from terms of the
 type $L_x\sigma_{<2}(D^2)$.\\
When $k_i=2$, i.e., terms  $$(L_x+
\Delta_x)^2\sigma_{<2}(D^2)=  L_x^2 \sigma_{<2}(D^2) + 2L_x\, \Delta_x\sigma_{<2}(D^2)+
\Delta_x^2\sigma_{<2}(D^2)$$ only contribute by $L_x^2\sigma_{<2}(D^2)$, which introduces terms of the type
\begin{eqnarray}\label{eq:Lx2} L_x^2 \left(\Gamma_{ij}^k\Gamma_{im}^n\sigma_{kj}\sigma_{nm}\right)&=&-4\,
\partial_a\Gamma_{ij}^k\partial_b
\Gamma_{im}^n \,\sigma_{kj}\sigma_{nm}\,\xi_a\xi_b\nonumber\\
&=& 
  -  \, R_{jkia} R_{mnib}\, \sigma_{jk} \sigma_{mn}\, \xi_a\,\xi_b.
\end{eqnarray}
 To sum up we only  have
contributions from  $\Delta_x\sigma_{<2}(D^2)$ and  $L_x^2\sigma_{<2}(D^2)$ via
products $$\left(L_x^2\sigma_{<2}(D^2)\right)^s \, (
\Delta_x\sigma_{<2}(D^2))^t\quad{\rm with}\quad  p=s+t \quad{\rm  and}\quad  \vert k\vert= 2s+t.$$   
Since the  residue picks the $-n$-th power in $\vert \xi\vert$, we have
 $$2s-2(\vert k\vert +q)=-n \Longrightarrow 2s+2t=\frac{n}{2}\Longrightarrow
 q=\frac{n}{4}.$$ This is confirmed by  counting the Clifford
 coefficients since 
 (\ref{eq:strgamma})   implies that $2q=2s+2t= \frac{n}{2}$. \\
 Combining this with (\ref{eq:Deltax}) and
(\ref{eq:Lx2}) leads to the statement of the proposition.
\endsquare
\begin{ex} When $n=4$,  in which case $q=1$, we have $s+t=1$ so that we need to consider
  two types of terms: $\Delta_x\sigma_{<2}(D^2)$ and $L^2_x(\sigma_{<2}(D^2))$.
   \\  
  Proposition \ref{prop:viaNCT} combined with Proposition
  \ref{prop:srespartialGamma} yields 
\begin{eqnarray*}&& {\rm sres}_x(\log(D^2))\\
&=& -\frac{1}{2}\,
{\rm sres}_x\left(\Delta_x(\sigma_{<2}(D^2))\,  \vert \xi\vert^{-4}\right)\quad
(s=0, t=1, k_1=1)\\
&+&  
\frac{{\rm sres}_x\left(L_x^2(\sigma_{<2}(D^2))\,  \vert \xi\vert^{-6}\right)}{2\times 3}\quad
(s=1, t=0,  k_1=2)\\
&=&{\rm sres}_x\left(  \vert
    \xi\vert^{-4}\,\partial_a \Gamma_{ij}^k
\partial_a\Gamma_{im}^n \,\sigma_{kj}\sigma_{nm}\right)\quad {\rm by}\quad (\ref{eq:Deltax})\\
&-&
\frac{4}{3} \, {\rm sres}_x\left( \vert
    \xi\vert^{-6}\,\partial_a \Gamma_{ij}^k
\partial_b\Gamma_{im}^n \,\sigma_{kj}\sigma_{nm} \,\xi_a \xi_b \right)\quad {\rm by}\quad
  (\ref{eq:Lx2})\\
&=&\left(1-\frac{1}{3}\right)\,{\rm sres}_x\left(  \vert
    \xi\vert^{-4}\,\partial_a \Gamma_{ij}^k
\partial_a\Gamma_{im}^n \,\sigma_{kj}\sigma_{nm}\right)\quad {\rm by}\quad (\ref{eq:srespartialGamma4bis})\\
&=& \frac{1}{48\, \pi^2}\,{\rm tr}\left( R\wedge R\right) \quad{\rm by}
\quad  (\ref{eq:srespartialGamma4}).\end{eqnarray*}
Once integrated over the manifold $M$ this yields back the well known formula: 
\begin{equation}\label{eq:AS4dim} {\rm ind}(D)= \int_M \hat A\Longrightarrow
{\rm sres}(\log (D^2))=-2\int_M \hat A= \frac{1}{48 \, \,\pi^2}
{\rm tr}\left( R\wedge R\right), \end{equation}
since  $\hat A= 1-\frac{1}{24 \, (2\,\pi)^2}{\rm tr}\left( R\wedge R\right)  
+\cdots $.  
\end{ex} 
\section*{Appendix: Complex powers and logarithms} 
An operator $A
\in \Cl(M,E)$ has principal angle $\theta$ if for every $(x, \xi)\in
T^*M-\{0\}$, the leading symbol $\left(\sigma_A(x, \xi)\right)^L$
has no eigenvalue on the ray  $L_\theta=\{re^{i\theta}, r\geq 0\}$;
in that case $A$ is
elliptic.
\begin{defn} We call an  operator $A\in \Cl(M, E)$
 admissible with spectral cut $\theta$
if $A$ has principal angle $\theta$ and the spectrum of $A$ does
 not meet $  {\rm L}_{\theta}=\{re^{i\theta}, r\geq  0\}$. In
particular such an operator is  invertible and elliptic. Since the
spectrum of $A$ does not meet $L_{\theta}$, $\theta$ is called an
Agmon angle of $A$.
\end{defn}
Let $A\in \Cl(M,E)$ be admissible with spectral cut $\theta$ and positive
order $a$. For ${\rm Re}(z)<0$, the complex power $A_{\theta}^z$ of $A$, first introduced by Seeley
\cite{Se},  is
defined by the Cauchy integral
\begin{equation} \label{eq:complexpowerCauchy}A_\theta^{z}=\frac{i}{2\pi} \int_{\Gamma_{r,\theta}} \lambda_\theta^{z}
(A-\lambda)^{-1}\, d\lambda,
\end{equation}
 where $\lambda_\theta^z= \vert
\lambda\vert^z e^{iz {(\rm arg}\lambda)}$ with  $\theta\leq {\rm
arg}\lambda<\theta+2\pi$. In particular,
 for $z=0$, we have $ A_{\theta}^0 = I$.
\\
Here
\begin{equation}\label{eq:contourGamma}\Gamma_{r,\theta}=\Gamma_{r,\theta}^1\cup\Gamma_{r,\theta}^2\cup\Gamma_{r,\theta}^3\end{equation}
where
 $$\Gamma_{r,\theta}^1=\{ \rho \,e^{i\theta}, \infty > \rho\geq r\}$$
 $$\Gamma_{r,\theta}^2=\{ \rho \,e^{i(\theta-2\pi)}, \infty > \rho\geq r\} $$
 $$\Gamma_{r,\theta}^3= \{r\, e^{it},\theta-2\pi\leq t \leq\theta\},$$
 is a contour along the ray $L_\theta$ around the non zero spectrum
of $A$. Here $r$ is any small positive real  number such that
$\Gamma_{r,\theta}\cap Sp(A) =\emptyset.$
\\ The definition of complex powers can be extended to the whole
complex plane by setting $A_{\theta}^z:=A^kA_{\theta}^{z-k}$ for
$k\in \N$ and ${\rm Re}(z)<k$; this definition is independent of the
choice of $k$ in $\N$ and preserves the usual properties, i.e.
$A_{\theta}^{z_1}A_{\theta}^{z_2}=A_{\theta}^{z_1+z_2}, $
$A_{\theta}^k=A^k, {\rm for}\quad k\in\Z.$ \\
Given  an admissible operator $A$ in $\Cl(M, E)$ with zero order  and   spectral cut
$\theta$, its complex powers give rise to a holomorphic map $z\mapsto
A_\theta^z$ on the complex plane with values in ${\cal B}(H^s(M, E))$ for any
real number $s$,
where $H^s(M, E)$ stands for the  $H^s$-closure of the
space $\Ci(M, E)$ of smooth sections of $E$ (se e.g. \cite{G}).  The
logarithm of $A$ is the bounded
operator  on $H^s(M, E)$  defined in terms of the derivative at
$z=0$
of this complex power:
\begin{eqnarray*}
\log_\theta A &:= & \left(\partial_z {A_\theta^{z}}\right)_{\vert_{z=0}}\\
&=& \frac{i}{2\pi}
\left(\partial_z\int_{\Gamma_{r,\theta}} \lambda_\theta^{z}\,(A-\lambda I)^{-z}\,
d\lambda\right)_{\vert_{z=0}}\\
&=& \frac{i}{2\pi}
\int_{\Gamma_{r,\theta}} \log_\theta \lambda\,(A-\lambda I)^{-z}\,
d\lambda
\end{eqnarray*}
with the notation of (\ref{eq:complexpowerCauchy}).\\
The notion of logarithm extends to an admissible operator $A$ with positive  order $a$ and   spectral cut
$\theta$ in the following way. For any positive $\e$, the map  $z\mapsto A_\theta^{z-\e}$
of order
$a(z-\e)$ defines a holomorphic  function on the half plane ${\rm Re}(z)<\e$ with
values  in ${\cal B}\left(H^s(M, E)\right)$ for any real number $s$.  Thus we can  set
\begin{equation}\label{eq:logCauchy} \log_\theta A =A_\theta^{\e}\,
 \left(\partial_z \left( A_\theta^z-\e\right)\right)_{\vert_{z=0}}= A_\theta^\e\,\left(\partial_z \left( \frac{i}{2\pi} \int_{\Gamma_{r,\theta}} \lambda_\theta^{z-\e}\,
(A-\lambda)^{-1}\, d\lambda\right)\right)_{\vert_{z=0}}.
 \end{equation}
 for any positive $\e$ the operator $ \log_\theta A\,  A^{-\e}=
A^{-\e}\, \log_\theta A$ lies in ${\cal B}\left(H^s(M, E)\right)$ for any
real number $s$. It follows  that $\log_\theta A $, which  is clearly independent of the choice of
  $\e>0$, defines a bounded
linear operator from $H^s(M, E)$ to $H^{s-\e}(M, E)$ for any positive $\e$.
WE have
$$\sigma(\log_\theta A)= a\, \log \vert \xi\vert + \sigma_0(\log_\theta A)$$
where $\sigma_0(\log_\theta A)$ is a classical symbol whose asymptotic
expansion $$\sigma_0(\log_\theta A)\sim\sum_{j=0}^\infty \sigma_{a-j,
  0}(\log_\theta A)$$   has homogeneous components of the form
\begin{equation}\label{eq:sigmaj0log}
\sigma_{-j, 0}(\log_\theta A)(x, \xi)=\vert \xi\vert^{-j}\,  \partial_z
\left(\sigma(A_\theta^z)_{\alpha(z)-j}(x,\frac{\xi}{\vert
  \xi \vert})\right)_{\vert_{z=0}}.
\end{equation}
\vfill\eject \noindent
\bibliographystyle{plain}

\end{document}